\documentclass[12pt]{article}
\usepackage{amssymb, amsmath, amsthm, graphics}
\newcommand{\RR}{{\mathbb R}}
\newcommand{\QQ}{{\mathbb Q}}

\newcommand{\ZZ}{{\mathbb Z}}
\newcommand{\FF}{{\mathbb F}}

\newtheorem{thm}{Theorem}[section]
\newtheorem{lem}[thm]{Lemma}
\newtheorem{prop}[thm]{Proposition}

\theoremstyle{definition}
\newtheorem*{rem}{Remark}
\newtheorem*{defn}{Definition}

\begin{document}

\title{\textbf{On embedding all $n$-manifolds \\into a single $(n+1)$-manifold}}
\author{Fan Ding, Shicheng Wang and Jiangang Yao}
\maketitle

\abstract{For each composite number $n\ne 2^k$, there does not
exist a single connected closed $(n+1)$-manifold such that any
smooth, simply-connected, closed $n$-manifold can be topologically
flat embedded into it.  There is a single connected closed
$5$-manifold $W$ such that any simply-connected, $4$-manifold $M$
can be topologically flat embedded into $W$ if $M$ is either
closed and indefinite, or compact and with non-empty boundary.}

\section{Introduction and some prerequisites}
\label{section:intro}

The celebrated Whitney embedding theorem states that any smooth
$n$-dimensional manifold can be embedded smoothly into the
Euclidean space $\RR^{2n}$, or equivalently, any smooth
$n$-dimensional manifold can be embedded into the sphere $S^{2n}$,
the simplest closed $2n$-manifold. If the target space is allowed
to be other closed manifolds, it is natural to wonder the
following general problem:

\textit{Find the smallest nonnegative integer $e_n$, such that any
$n$-dimensional connected, closed manifold can be (topologically
flat or smoothly) embedded into a single connected, closed
manifold of dimension $n+e_n$.}

Clearly $0 \le e_n \le n$.

{\bf Examples.} The only known $e_i$'s are for $i=0, 1, 2, 3$.

(i) $e_0=0$ since a connected 0-manifold is a point;

(ii) $e_1=0$ since a connected, closed 1-manifold is homeomorphic
to $S^1$;

(iii)  $e_2=1$ since any orientable, connected, closed surface
embeds into $S^3$, and any non-orientable, connected, closed
surface can be presented as a connected sum of an orientable,
connected, closed surface and $\mathbb{RP}^2$ or
$\mathbb{RP}^2\#\mathbb{RP}^2$, hence each connected closed
surface can be embedded into $S^3\# \mathbb{RP}^3\#
\mathbb{RP}^3\cong \mathbb{RP}^3\#\mathbb{RP}^3$.

(iv) $e_3=2$ since each closed 3-manifold embeds into $S^5$, due
to Hirsch in the orientable case ([Hi], also [Ki]), and due to
Rohlin and Wall independently in the non-orientable case ([Ro]
[Wa], see [WZ] for a possibly more elementary proof, in the sense
that only Dehn surgery is involved, of this Hirsch-Rohlin-Wall's
Theorem); and Kawauchi showed that there does not exist a single
oriented, connected, closed $4$-manifold such that any connected,
closed $3$-manifold can be topologically flat embedded into it by
constructing signature invariants from infinite cyclic coverings
[Ka], and then Shiomi showed that Kawauchi's result is still true
if the target 4-manifold is allowed to be non-orientable [Shi].

\medskip
It is reasonable to guess that in general $e_n>1$. The main result
of the present paper claims that this is the case for ``most" $n$,
and indeed we prove a stronger version for those $n$'s.

\textbf{Theorem~\ref{thm:embed}.} \textit{If $n$ is a composite
number and is not a power of $2$, then there does not exist a
single connected, closed $(n+1)$-manifold $W$, such that any
smooth, simply-connected, closed $n$-manifold $M$ can be
topologically flat embedded into $W$.}

We say that $M$ can be \textit{topologically flat} embedded into
$W$ if $M\times [0,1]$ can be embedded into $W$.

Besides [Ka] and [Shi], recently there are several papers studying
the embedding of a 3-manifold $M$ into simply-connected
4-manifolds other than $S^4$, say [EL], [Fa1], [Fa2], and indeed
the targets in those papers are connected sums of $\mathbb{CP}^2$.
This relates to finding uncountably many smooth structures on the
4-manifold $M\times \RR$ (see [Di], [Fa1] and [Fa2]). Fang's
observation (cf. [Fa2]) that there exists an obstruction to the
embedding of 3-manifolds into a simply connected definite
4-manifold was the first inspiration of our proof for
Theorem~\ref{thm:embed}.

\medskip
We expect that Theorem~\ref{thm:embed} still holds for all
remaining positive integers $n\geq 5$. When $n=3$, the only
simply-connected, closed 3-manifold is $S^3$ if Poincar\'e
Conjecture is valid. When $n=4$, the situation is still not clear,
and we have the following

\textbf{Theorem~\ref{thm:dim4case}.} \textit{(a) There exists an
oriented, connected, closed $5$-manifold $W$, such that for any
simply-connected, closed $4$-manifold $M$ with indefinite
intersection form, $M$ can be topologically flat embedded into
$W$.}

\textit{(b) Any simply-connected, compact topological $4$-manifold
$M$ with $\partial M \not = \emptyset$ admits a topologically flat
embedding into $S^2 \tilde\times S^3$.}

\vspace{-5pt}
\begin{rem} There are indefinite 4-manifolds with arbitrarily large
signatures. Theorem~\ref{thm:dim4case}~(a) tells us that the
signature of 4-manifolds cannot be an obstruction for our
codimension 1 embeddings. Note that any oriented, connected,
closed 4-manifold with non-zero signature does not bound a compact
orientable 5-manifold. It follows that the first Betti number of
$W$ in (a) must be positive and that the condition $\partial M\ne
\emptyset$ in (b) can not be removed.
\end{rem}
\vspace{-5pt}

Compared with Theorem~\ref{thm:dim4case}~(b), it is worthy noting
that there does not exist a single oriented, connected, closed
$4$-manifold such that any compact punctured $3$-manifold can be
topologically flat embedded into it (cf. [Ka]).

\medskip

\textbf{The organization of the paper and an outline of the
proofs.}

In Section~\ref{section:obstr}, we first reduce the proof of
Theorem~\ref{thm:embed} to the oriented case. Then for any fixed
(n+1)-manifold $W$, Proposition~\ref{prop:obstr} claims that for
each topologically flat codimension 1 embedding $M^n\to W$ and
each factorization $n=pq$, there is a homological obstruction in
terms of $H^{p+1}(W;\RR)$, $H^p(M;\RR)$, the $q$-multiple cup
product on $H^p(M;\RR)$ and so on. Since $n$ is not a prime number
and is not a power of 2, we can assume that the factorization
$n=pq$ has been chosen so that $p\ge 2$ and $q\ge 3$, and $q$ is
odd.

Let $(\wedge^q\RR^m)^*$ (resp. $(\vee^q\RR^m)^*$) be the vector
space of all skew-symmetric (resp. symmetric) $q$-multilinear
functions on the $m$-dimensional real vector space $\RR^m$.
Motivated by the homological obstruction in
Proposition~\ref{prop:obstr}, in Section~\ref{section:special} we
define when a  function in $(\wedge^q\RR^m)^*$ (resp.
$(\vee^q\RR^m)^*$) is special so that the $q$-multiple cup product
on $H^p(M;\RR)$ in Proposition~\ref{prop:obstr} is a special
function. Then Proposition~\ref{prop:special} claims that under
certain circumstances, such special functions are contained in a
proper closed subset of $(\wedge^q\RR^m)^*$ (resp.
$(\vee^q\RR^m)^*$). The conditions that $q\ge 3$, and $q$ is odd
are used in proving Proposition~\ref{prop:special}.

In Section~\ref{section:main}, we first choose a suitable rational
non-special function provided by Proposition~\ref{prop:special},
and then Proposition~\ref{prop:extend} claims that there is a
commutative graded algebra $A=\bigoplus_{i=0}^{q} A_{ip}$ over
$\QQ$ satisfying Poincar\'e duality so that the chosen non-special
function is the $q$-multiple product on $A_p$. Next we verify that
$A$ can be realized as the rational cohomology ring of a
simply-connected closed $n$-manifold $M$ by invoking a celebrated
result of Sullivan (the condition $p\ge 2$ is used here). Hence
the $q$-multiple cup product on $H^p(M;\RR)$ is a non-special
function, and therefore $M$ cannot be topologically flat embedded
into $W$ by Proposition~\ref{prop:obstr}, which finishes the proof
of Theorem~\ref{thm:embed}.

Section~\ref{section:dim4case} is devoted to the case $n=4$, and
we get Theorem~\ref{thm:dim4case} by explicit constructions and
the classification theorem for simply-connected, closed
$4$-manifolds.

\medskip
In the remaining of this section, We recall and fix needed
notations and conventions in tensor algebra  based on [Sh].

{\bf Prerequisites on tensor algebra.} Let $V$ be a vector space
of finite dimension $m$ over a field $\FF$ of characteristic zero.
The space $\otimes ^p V$, consisting of contravariant tensors of
order $p$, is defined for $p=2,3,\ldots$ to be the tensor product
$V\otimes V\otimes\ldots\otimes V$ of $p$ copies of $V$, and for
$p=0,1$ to be $\FF ,V$. We may identify
$(\otimes^pV)\otimes(\otimes^qV)$ with $\otimes^{p+q}V$. Thus
given contravariant tensors $x,y$ of orders $p,q$ we may form
their product $x\otimes y$, a contravariant tensor of order $p+q$.
Observe that the product of decomposable elements is given by
$(x_1\otimes\ldots\otimes x_p)\otimes (y_1\otimes\ldots\otimes
y_q)=x_1\otimes\ldots\otimes x_p\otimes y_1\otimes\ldots\otimes
y_q$. Form the direct sum
$$\otimes V=\oplus_p(\otimes ^p V)=\FF\oplus V\oplus
(\otimes^2V)\oplus\ldots$$ of all the spaces
$\otimes^pV,p=0,1,2,\ldots$. If we extend the product operation
$\otimes$ in the obvious bilinear way to general elements of the
vector space $\otimes V$, we make $\otimes V$ into an associative
algebra over $\FF$, called the \textit{tensor algebra} of the
space $V$.

Let $N^p(V)$ denote the subspace of $\otimes ^pV$ spanned by those
tensors of the form $v_1\otimes\ldots\otimes v_p$, where
$v_i=v_{i+1}$ for some $i=1,2,\ldots ,p-1$. For $p=0$ or $1$,
$N^p(V)$ is defined to be $0$. The $p$th \textit{exterior power}
of $V$ is defined to be the quotient space
$\wedge^pV=\otimes^pV/N^p(V)$. Observe that $\wedge^0V=\FF$,
$\wedge^1V=V$. The canonical projection $\pi_p:\otimes^pV\to
\wedge^pV$ is defined by $\pi_p(x)=x+N^p(V)$ for $x\in\otimes^pV$.
For $v_1,\ldots,v_p\in V$, we will write $\pi_p(v_1\otimes
\ldots\otimes v_p)$ as $v_1\wedge \ldots \wedge v_p$. There is a
product operation
$\wedge:(\wedge^pV)\times(\wedge^qV)\to\wedge^{p+q}V$ defined as
follows. If $x\in\otimes^pV,y\in\otimes^qV$, then the
\textit{exterior product} of $\pi_p(x)$ and $\pi_q(y)$ is defined
by $\pi_p(x)\wedge \pi_q(y)=\pi_{p+q}(x\otimes y)$. Observe that
the exterior product of decomposable elements is given by
$(x_1\wedge\ldots\wedge x_p)\wedge (y_1\wedge\ldots\wedge
y_q)=x_1\wedge\ldots\wedge x_p\wedge y_1\wedge\ldots\wedge y_q$.
If we extend $\wedge$ in the obvious bilinear way to general
elements of the vector space
$$\wedge V=\oplus_p(\wedge^p V)=\FF\oplus
V\oplus\wedge^2V\oplus\ldots\oplus\wedge^mV,$$ we make $\wedge V$
into an algebra over $\FF$, called the \textit{exterior algebra}
on $V$.

Let $S_p$ denote the symmetric group on the $p$ symbols $\{
1,2,\ldots,p\}$. Each permutation $\sigma\in S_p$ gives rise to a
linear operator, still denoted by $\sigma$, on $\otimes^pV$
defined by its effect $$\sigma(v_1\otimes\ldots\otimes
v_p)=v_{\sigma^{-1}(1)}\otimes\ldots\otimes v_{\sigma^{-1}(p)}$$
upon the decomposable elements of $\otimes^pV$. Let $K^p(V)$
denote the subspace of $\otimes ^pV$ spanned by those tensors of
the form $u-\sigma u$ for some $u\in\otimes^pV$ and  some
$\sigma\in S_p$. If $p=0$ or $1$ we define $K^p(V)$ to be $0$. The
$p$th \textit{symmetric power} of $V$ is defined to be the
quotient space $\vee^pV=\otimes^pV/K^p(V)$. Observe that
$\vee^0V=\FF$, $\vee^1V=V$. The canonical projection
$\psi_p:\otimes^pV\to \vee^pV$ is defined by $\psi_p(x)=x+K^p(V)$
for $x\in\otimes^pV$. Then exactly parallel to the description in
the preceding paragraph (simply changing $\pi_p$ to $\psi_p$ and
$\wedge$ to $\vee$), we see that there is a product operation
$\vee :(\vee^pV)\times(\vee^qV)\to\vee^{p+q}V$ which can be
extended in the obvious bilinear way to general elements of the
vector space
$$\vee V=\oplus_p(\vee^p V)=\FF\oplus
V\oplus\vee^2V\oplus\ldots$$ to make $\vee V$ into an
infinite-dimensional algebra over $\FF$, called the
\textit{symmetric algebra} on $V$.

Let $\times^p V$ denote the Cartesian product $V\times
V\times\ldots\times V$ of $p$ copies of $V$. Each multilinear
function $F:\times^p V\to \FF$ determines a unique linear function
$L:\otimes^pV\to \FF$ such that for all $v_1,\ldots,v_p\in V$,
$L(v_1\otimes\ldots\otimes v_p)=F(v_1,\ldots,v_p)$. The
correspondence $F\leftrightarrow L$ establishes a natural
isomorphism between the vector space $L(\times^pV;\FF)$ of all
multilinear functions and the dual space $(\otimes^pV)^*$ of
$\otimes^pV$. For any skew-symmetric (resp. symmetric) multilinear
function $F:\times^p V\to \FF$, there exists a unique linear
function $L_s:\wedge^pV$ (resp. $\vee^pV$) $\to \FF$ such that for
all $v_1,\ldots,v_p\in V$, $L_s(v_1\wedge\ldots\wedge
v_p)=F(v_1,\ldots,v_p)$ (resp. $L_s(v_1\vee\ldots\vee
v_p)=F(v_1,\ldots,v_p)$). The correspondence $F\leftrightarrow
L_s$ establishes a natural isomorphism between the vector space
${\rm Sk}(\times^pV ;\FF)$ (resp. ${\rm Sym}(\times^pV;\FF)$) of
all skew-symmetric (resp. symmetric) multilinear functions and the
dual space $(\wedge^pV)^*$ (resp. $(\vee^pV)^*$) of $\wedge^pV$
(resp. $\vee^pV$). Thus we may consider a skew-symmetric (resp.
symmetric) multilinear function $F:\times^p V\to \FF$ as an
element of $(\wedge^pV)^*$ (resp. $(\vee^pV)^*$), and vice versa.

Let $\{ e_1,\ldots,e_m\}$ be the standard basis for the real
vector space $\RR^m$, i.e. $e_i=(0,\ldots,0,1,0,\ldots,0)$ has
components 0 except for its i-th component, which is equal to 1.
Note that $(\wedge^q\RR^m)^*$ (resp. $(\vee^q\RR^m)^*$) has
dimension
$\begin{pmatrix} m \\ q  \end{pmatrix}$ (resp. $\begin{pmatrix} m+q-1 \\
q\end{pmatrix}$). Each element $F\in (\wedge^q\RR^m)^*$ (resp.
$(\vee^q\RR^m)^*$) is determined by $\begin{pmatrix} m \\ q
\end{pmatrix}$ (resp. $\begin{pmatrix} m+q-1 \\
q\end{pmatrix}$) real numbers $(F_{i_1\ldots i_q})$, where
$F_{i_1\ldots i_q}=F(e_{i_1},\ldots,e_{i_q})$ for $1 \le i_1<
\ldots< i_q\le m$ (resp. $1\le i_1\le \ldots\le i_q\le m$). We
call $F_{i_1\ldots i_q}$ the \textit{$(i_1,\ldots,i_q)$ component}
of $F$.

\section{Reduction to the oriented case and homological obstructions for the embedding}
\label{section:obstr}

If any smooth, simply-connected, closed $n$-manifold $M$ can be
topologically flat embedded into a non-orientable $(n+1)$-manifold
$W$, then $M$ can also be topologically flat embedded into the
orientable double cover of $W$ because $M$ is simply-connected. It
follows that Theorem~\ref{thm:embed} is equivalent to the
following

\textbf{Theorem~\ref{thm:embed}*.} \textit{If $n$ is a composite
number and is not a power of $2$, then there does not exist a
single oriented, connected, closed $(n+1)$-manifold $W$, such that
any smooth, simply-connected, closed $n$-manifold $M$ can be
topologically flat embedded into $W$.}

From now on all manifolds appearing in this paper are considered
to be oriented and connected. $\beta_i(X)$ denotes the $i$-th
Betti number of a compact manifold $X$.

The following proposition sets up a property on the cohomology
ring of $n$-manifold $M$ if it can be embedded into the
$(n+1)$-manifold $W$.

\begin{prop}
\label{prop:obstr} Let $M$ and $W$ be closed, connected, oriented
$n$- and $(n+1)$-dimensional  manifolds respectively.

If $M$ topologically flat embeds in $W$, then for any integer
factorization $n=pq$, where $p,q>0$, there exists a subspace $V$
of $H^p(M;\RR)$ and a linear transformation $\varphi: V \to
H^p(M;\RR)$ such that

(i) $\varphi$ has  no fixed non-zero vectors, i.e. for any
non-zero $x\in V$, $\varphi (x)\ne x$,

(ii) $\dim V \geq \frac{1}{2}(\beta_p(M)-\beta_{p+1}(W))$,

(iii) for any $x_1, \ldots, x_q \in V$, $x_1 \cup \ldots \cup x_q
=  \varphi(x_1) \cup \ldots \cup \varphi(x_q).$
\end{prop}

\begin{proof} All homology and cohomology groups in this proof are with
coefficients in $\RR$. We divide the proof into two cases.

\underline{Case 1}. Suppose $W \backslash M$ is not connected.
Denote the closures of the two components of $W \backslash M$ by
$W_1$ and $W_2$. By Mayer-Vietoris sequence, we know that
$$H^p(W_1) \oplus H^p(W_2) \stackrel{f}\to H^p(M) \stackrel{\delta^*} \to H^{p+1}
(W)$$ is exact. Here the first map is given by $f(w_1,w_2)=
i_1^*(w_1)-i_2^*(w_2)$ for $w_1\in H^p(W_1),w_2\in H^p(W_2)$,
where $i_1$ and $i_2$ are inclusions $M \hookrightarrow W_1$ and
$M \hookrightarrow W_2$ respectively. Clearly
$$\dim {\rm im} (f) =\dim \ker (\delta^*) \geq  \dim H^p(M) - \dim
H^{p+1}(W) =\beta_p(M)-\beta_{p+1}(W).$$ Since $\dim {\rm im}
(i_1^*) + \dim {\rm im}(i_2^*) \geq \dim {\rm im} (f)$, we may set
$V= {\rm im}(i_1^*)$ and assume that $\dim V \geq
\frac{1}{2}(\beta_p(M)-\beta_{p+1}(W))$.

For any $x_1, \ldots, x_q \in V$, there exist $y_1, \ldots, y_q
\in W_1$ such that $x_i=i_1^*(y_i), 1\leq i \leq q$. Since
$M=\partial W_1$, we have $i_{1*}[M]=0 \in H_n(W_1)$, where $[M]$
denotes the fundamental class of $M$. So
$$ \begin{array}{lll}
<x_1 \cup \ldots \cup x_q, [M]> & = & < i_1^*(y_1 \cup \ldots \cup y_q), [M]> \\
  & = & < y_1\cup \ldots \cup y_q, i_{1*}[M]>=0.
\end{array}$$
Therefore for any $x_1, \ldots, x_q \in V$, $x_1 \cup \ldots \cup
x_q=0.$

The result is thus proved by letting $\varphi$ to be the zero map.

\underline{Case 2}. Suppose $W \backslash M$ is connected. We
identify a neighborhood $W_1$ of $M$ in $W$ with $M \times I$ and
let $W_2=W \backslash {\rm int}(W_1)$, where ${\rm int}(W_1)$
denotes the interior of $W_1$. Clearly, $W_1 \cap W_2= (M\times \{
0 \}) \cup (M \times \{ 1 \})$.

Still by Mayer-Vietoris sequence, we get the exact sequence
$$H^p(W_1) \oplus H^p(W_2) \to H^p(W_1 \cap W_2)  \to H^{p+1}
(W)$$ and it is reduced to
$$H^p(M) \oplus H^p(W_2) \stackrel{f} \to H^p(M)\oplus H^p(M)
\stackrel{\delta^*} \to H^{p+1}(W).$$ Now $f$ can be written as
$$f(x,w)=(x-i_0^*(w), x-i_1^*(w)),$$
where $i_0$ and $i_1$ are inclusions $M \times \{ 0 \}
\hookrightarrow W_2$ and $M \times \{ 1 \} \hookrightarrow W_2$
respectively. By the exactness,
$$\dim {\rm im} (f) =\dim \ker (\delta^*) \geq  \dim (H^p(M)\oplus H^p(M)) - \dim
H^{p+1}(W) = 2\beta_p(M)-\beta_{p+1}(W).
$$

Let $ \Delta = \{(x,x)|x\in H^p(M)\} \subset H^p(M)\oplus H^p(M)$
and $U=\{ (i_0^*(w), i_1^*(w)) | w \in H^p(W_2) \}$, then $\Delta+
U={\rm im}(f)$ and therefore
$$\dim (\Delta+ U)=\dim {\rm
im}(f)\ge 2\beta_p(M)-\beta_{p+1}(W).$$ Since $\dim \Delta =
\beta_p(M)$, there is a subspace $S$ of $U$ with $\dim S \geq
\beta_p(M)-\beta_{p+1}(W)$ such that $S\cap \Delta=\{(0,0)\}$. Let
$\pi_1,\pi_2:H^p(M)\oplus H^p(M)\to H^p(M)$ be the projections of
$H^p(M)\oplus H^p(M)$ onto its first and second factors,
respectively. Since $\dim \pi_1(S)+ \dim \pi_2(S) \geq \dim S$,
without loss of generality, we may assume that
$$\dim \pi_1(S) \geq \frac{1}{2}\dim S \ge \frac{1}{2}
(\beta_p(M) - \beta_{p+1}(W)).$$

It follows from $\partial W_2= -( M \times \{0\}) \cup (M \times
\{1\})$ that $i_{0*}[M]=i_{1*}[M]$. Thus for any $w_1, \ldots, w_q
\in H^p(W_2) $, $< w_1\cup \ldots \cup w_q,
i_{0*}[M]-i_{1*}[M]>=0$, and this implies
$$ < i^*_0(w_1) \cup \ldots \cup i^*_0(w_q), [M]> =< i^*_1(w_1) \cup \ldots \cup i^*_1(w_q),
[M]>.$$ In particular, if $(u_1,v_1),\ldots,(u_q,v_q)\in S$, then
\begin{equation} u_1\cup\ldots\cup u_q=v_1\cup\ldots\cup v_q.
\label{eq:cupidentity} \end{equation}

Let $\{ y_1, \ldots, y_r \}$ be a basis for $\pi_1(S)$. We can
find $z_1, \ldots, z_r\in H^p(M)$ such that $(y_j,z_j)\in S$,
$1\le j\le r$. Then there is a unique linear transformation
$$\varphi: \pi_1(S) \to H^p(M)$$ such that
$\varphi(y_j)=z_j, 1\leq j \leq r$.

One can write any $y \in \pi_1(S)$ as $y = \sum_{j} c_j y_j$, then
$\varphi(y)=\sum_{j} c_jz_j $. Thus $(y,\varphi(y))\in S$. Since
$S\cap\Delta=\{ (0,0)\}$, if $y\ne 0$, then $\varphi(y)\ne y$. Now
set $V=\pi_1(S)$ and one gets the desired result from
equation~(\ref{eq:cupidentity}).
\end{proof}

\section{Special multilinear functions}
\label{section:special}

Proposition~\ref{prop:obstr} motivates the following definition.

\begin{defn} Let $F$ be an element of $(\wedge^q\RR^m)^*$ (resp.
$(\vee^q\RR^m)^*$). We say that $F$ is \textit{special} if there
exist a subspace $U$ of $\RR^m$ with $\dim U\ge \frac{m}{3}$ and a
linear map $\varphi: U \to \RR^m$ with no fixed non-zero vectors
(i.e. for any non-zero $x\in U$, $\varphi(x)\ne x$) such that for
all $x_1,\ldots,x_q\in U$,
$$F(x_1, \ldots, x_q)= F( \varphi(x_1), \ldots,\varphi(x_q)).$$
In the above equation, we consider $F$ as a skew-symmetric (resp.
symmetric) multilinear function (cf. Section~\ref{section:intro}).
\end{defn}

In this section, we want to show

\begin{prop}
\label{prop:special}
Suppose $q \geq 3$ is an odd integer.

(a) If $m$ is sufficiently large, then there exists a proper
closed subset $X_m$ of $(\wedge^q\RR^m)^*$ such that if $F \in
(\wedge^q\RR^m)^*$ is special, then $F\in X_m$.

(b) The result in (a) is still true if we substitute $(\wedge
^q\RR^m)^*$ by $(\vee ^q\RR^m)^*$.
\end{prop}

\begin{rem} If $q\ge 2$ is even, then each $F\in (\wedge^q\RR^m)^*$
(resp. $F\in (\vee^q\RR^m)^*$) is special. We may simply take
$U=\RR^m$ and $\varphi:\RR^m\to\RR^m$ to be the map sending $x$ to
$-x$ for all $x\in \RR^m$.
\end{rem}

Before proving Proposition~\ref{prop:special}, we shall discuss
some useful lemmas. The first lemma describes the real Jordan
canonical form.

\begin{lem}
\label{lem:Jordan}
Any real square matrix is similar to a block
diagonal matrix, where each block takes one of the following two
forms
$$J_k(a)=\begin{bmatrix} a & 1 & &  & \\ & a & 1 & &
 \\ & & \ddots& \ddots &  \\ & & & \ddots& 1  \\ & & &
 & a  \end{bmatrix},
C_l(a,b)=\begin{bmatrix} a & b &  1 & 0 & & & &  \\
-b & a & 0 & 1 & & & & \\
& & a & b & \ddots & & & \\
& & -b & a & & \ddots & & \\
& & & & \ddots & & 1 & 0  \\
& & & & & \ddots &  0 & 1 \\
& & & & & & a & b \\
& & & & & & -b & a
\end{bmatrix}. $$
The sizes of these two square matrices are $k$ and $2l$
respectively, where $k$ and $l$ are positive integers, $a$ and
$b\not =0$ are real numbers.
\end{lem}

\begin{proof} This is Theorem~3.4.5 in [HJ]. \end{proof}

Now we sort the blocks described in this lemma into 6 types for
later use.
$$\begin{array}{ccc} \textbf{Type 1:} \ J_1(1), & \textbf{Type 2:}  \ J_1(-1), & \textbf{Type 3:} \ J_1(a),
|a|>1, \\
\textbf{Tpye 4:} \ J_1(a), |a|<1, & \textbf{Type 5:} \ J_k(a),
k>1, & \textbf{Type 6:} \ C_l(a,b). \end{array}$$

\begin{lem}
\label{lem:subinj}
Let $V$ be a subspace of $\RR^m$ with $\dim
V=l$. Suppose that $\{ v_1, \ldots, v_l \}$ is a basis for $V$,
and that $ \{ v_1, \ldots, v_l, \ldots, v_m \}$ is a basis for
$\RR^m$.

Let $\varphi: V \to \RR^m$ be a linear map with no fixed non-zero
vectors such that the matrix of $\varphi$ with respect to the
above bases is $\begin{bmatrix} J
\\ *
\end{bmatrix}$, where $J$ is an $l \times l$ block diagonal matrix
and all the blocks are of Type $i_0$, $i_0=1,5$ or $6$. Then there
exists a subspace $V'$ of $V$ such that $\dim V' \geq
\frac{l}{3},\ \varphi(V') \cap V' =\{ 0 \}$ and $\varphi$ is
injective on $V'$.
\end{lem}

\begin{proof} First suppose that $J$ only contains Type 1 blocks. Let
$W={\rm span}\{v_{l+1}, \ldots, v_{m}\}$, then $\RR^m= V\oplus W$.
For any $v \in V$, since $J=I$, we have decomposition
$$\varphi(v)=v + w, v \in V, w \in W.$$
Thus $\varphi$ is injective. Suppose that $x \in \varphi(V) \cap
V$, then $x\in V$ and there exists $v\in V$ such that
$\varphi(v)=x$. Now the decomposition of $\varphi(v)$ under $V
\oplus W$ equals both $x+0$ and $v+w$, so $x=v$, $\varphi(v)=v$.
Since $\varphi$ only fixes the zero vector, we have $v=0$.
Therefore $\varphi(V) \cap V = \{ 0 \}$, and we may simply choose
$V'=V$.

Suppose now that $J$ only contains Type 5 or Type 6 blocks. For a
block in $J$, let $\{ w_1, \ldots, w_{k} \}$ be the subset of the
basis $\{ v_1, \ldots, v_l\}$ that corresponds to that block. Note
that $k\geq 2$. Set $W={\rm span}\{ w_1, \ldots, w_k\}$, $W_0={\rm
span} \{ w_j |$ $j \ {\rm even} \}$ and $W_1={\rm span}\{ w_j | $
$j \ {\rm odd} \}$. Clearly, $W=W_0 \oplus W_1$ and $\dim {W_0}
\geq \frac{k}{3}$. Denote the projection of $\RR^m$ onto $W_1$ by
$p_1$. (This projection is well-defined as there is a natural
decomposition $\RR^m=W\oplus W'$, where $W'={\rm span}(\{ v_1,
\ldots, v_m\}\backslash\{w_1, \ldots, w_{k}\})$.) If the block is
$C_\frac{k}{2}(a,b)$, then \begin{equation} p_1\varphi(w_{2i})=
bw_{2i-1},\ 1 \leq i \leq \frac{k}{2}, \label{eq:p1C}
\end{equation} and if the block is $J_k(a)$, then
\begin{equation} p_1\varphi(w_{2i})=w_{2i-1},\  1 \leq i \leq \left [\frac{k}{2}\right ],
\label{eq:p1J} \end{equation} where $[ \frac{k}{2}]$ is the
greatest integer $\le \frac{k}{2}$. Note that $b\not =0$. In any
case, $p_1\varphi(x)=0$ holds for $x\in W_0$ only if $x=0$. Thus
$\varphi$ is injective on $W_0$. By equations~(\ref{eq:p1C}) and
(\ref{eq:p1J}), $\varphi(W_0) \cap W_0 =\{ 0 \}$. Now we may
choose $V'$ to be the direct sum of these $W_0$'s of all blocks in
$J$.
\end{proof}

The next lemma needs to use the transcendence degree of a field
extension. A good reference is [La].

\begin{lem} \label{lem:algdep} Let $g_1(x_1, \ldots,
x_r),\ldots,g_s(x_1,\ldots,x_r)$ be rational functions in $r$
variables $x_1, \ldots, x_r$ with coefficients in $\RR$. If $s>r$,
then there exists a non-zero polynomial $P$ in $s$ variables such
that
$$P(g_1(x_1 \ldots, x_r), \ldots, g_s(x_1,
\ldots, x_r)) = 0.$$
\end{lem}

\begin{proof} Let $K$ be the field of fractions of the polynomial ring
$\RR[x_1, \ldots, x_r]$. Clearly, $\{ x_1, \ldots, x_r \}$ is a
transcendence base of $K$ over $\RR$. Since $s>r$, $\{ g_1,
\ldots, g_s\}$ is not algebraically independent over $\RR$ (cf.
[La], Chapter~VIII, Theorem~1.1). Hence there exists a non-zero
polynomial $P$ in $s$ variables with coefficients in $\RR$ such
that $ P(g_1, \cdots, g_s)=0.$ \end{proof}

\textbf{Proof of Proposition~\ref{prop:special}.}

(a) Suppose that $F_0\in (\wedge^q\RR^m)^*$ is special. Thus there
exist a subspace $U$ of $\RR^m$ with $\dim U =k\geq \frac{m}{3}$
and a linear map $\varphi: U \to \RR^m$ with no fixed non-zero
vectors such that for all $x_1,\ldots,x_q\in U$,
$$F_0(x_1, \ldots, x_q)= F_0( \varphi(x_1), \ldots,\varphi(x_q)).$$

We will use the above equation to get expressions for coordinates
of $F_0$ under standard basis for $(\wedge^q\RR^m)^*$. The number
of variables used will be less than $\dim (\wedge^q\RR^m)^*$ when
$m$ is large enough . Hence all special $F_0$ lie in the zero set
of a single non-zero polynomial, which is a proper closed subset
of $(\wedge^q\RR^m)^*$.

Choose vectors $u_{k+1}, \ldots, u_{m}$ such that
$$U\oplus {\rm span}\{u_{k+1}, \ldots, u_{m} \}=\RR^m.$$
Let $\{u_1, \ldots, u_{k}\}$ be a basis for $U$, then $\{u_1,
\ldots, u_{k}, u_{k+1}, \ldots, u_{m}\}$ is a basis for $\RR^m$.
Denote the matrix of $\varphi$ with respect to these two bases by
$\begin{bmatrix} J
\\ *
\end{bmatrix}$, where $J$ is a $k\times k$ square matrix.
A suitable choice of $\{ u_1, \ldots, u_k\}$ can make $J$ in real
Jordan canonical form as in Lemma~\ref{lem:Jordan}. The Jordan
blocks have 6 types. For some $i_0$, the total size of blocks of
Type $i_0$ in $J$ should be at least $\frac{k}{6} \geq
\frac{m}{18}$. Denote this total size by $l$. We change the order
of $u_1, \ldots, u_k$ so that in the block expression for $J$,
Type $i_0$ blocks appear first. From now on, vectors $u_1, \ldots,
u_k, \ldots, u_m$ are fixed.

\underline{Case 1}. Suppose that $i_0=1,5$ or $6$. Then by
Lemma~\ref{lem:subinj}, we may find a subspace $U'$ of ${\rm
span}\{u_1,\ldots,u_{l}\}$ with $k'={\rm dim} U'\ge \frac{l}{3}\ge
\frac{m}{54},\ \varphi(U')\cap U'=\{ 0\}$ and $\varphi$ is
injective on $U'$. Let $m_1$ denote $[\frac{m}{54}]$, the greatest
integer $\le \frac{m}{54}$. Let $\{v_1,\ldots,v_{k'}\}$ be a basis
of $U'$. Then $\{
v_1,\ldots,v_{k'},\varphi(v_1),\ldots,\varphi(v_{k'})\}$ is
linearly independent and there exists an isomorphism
$\tau_1:\RR^m\to \RR^m$ sending $v_{i}$ to $e_{i}$, $\varphi(v_i)$
to $e_{m-m_1+i}$ for $1 \leq i \leq m_1$. (cf. section 1 for
definition of $e_i$.)

Let $L(\RR^m,\RR^m)$ be the vector space consisting of all linear
maps from $\RR^m$ to $\RR^m$. $L(\RR^m,\RR^m)$ has dimension
$m^2$. Each element $\tau\in L(\RR^m,\RR^m)$ is determined by
$m^2$ numbers $(\tau_{ji})$, where $\tau(e_i)=\sum_{j=1}^m
\tau_{ji}e_j$ for $1\le i\le m$. For each $\tau\in
L(\RR^m,\RR^m)$, $\tau$ induces a linear map
$\tau^*:(\wedge^q\RR^m)^*\to (\wedge^q\RR^m)^*$ defined by
$$(\tau^*F)(x_1, \ldots, x_q)= F( \tau (x_1), \ldots, \tau
(x_q))$$ for $F \in (\wedge^q\RR^m)^*$ and $x_1, \ldots, x_q \in
\RR^m$.

Let $R=\{ (i_1, \ldots, i_q) | 1\leq i_1< \ldots< i_q \leq m \}$,
$R_1=\{ (i_1, \ldots, i_q) | 1\leq i_1< \ldots< i_q \leq m_1\}$
and $R_1^{\prime}=R\backslash R_1$.

Now consider $(\tau_1^{-1})^*F_0$. Note that for
$(i_1,\ldots,i_q)\in R_1$, we have
$$((\tau_1^{-1})^*F_0)_{i_1\ldots
i_q}=((\tau_1^{-1})^*F_0)(e_{i_1},\ldots,e_{i_q})$$
$$=F_0(v_{i_1},\ldots,v_{i_q})=F_0(\varphi(v_{i_1}),\ldots,\varphi(v_{i_q}))$$
$$=((\tau_1^{-1})^*F_0)(e_{m-m_1+i_1},\ldots,e_{m-m_1+i_q})$$
$$=((\tau_1^{-1})^*F_0)_{(m-m_1+i_1)\ldots
(m-m_1+i_q)}.$$

Let $T$ denote the subspace of $(\wedge^q\RR^m)^*$ consisting of
all $F\in (\wedge^q\RR^m)^*$ such that $F_{i_1\ldots
i_q}=F_{(m-m_1+i_1)\ldots (m-m_1+i_q)}$ for all
$(i_1,\ldots,i_q)\in R_1$. Note that $(\tau_1^{-1})^*F_0\in T$. As
each element $F\in T$ is determined by $\begin{pmatrix} m
\\ q  \end{pmatrix}-\begin{pmatrix} m_1 \\
q \end{pmatrix}$ numbers $(F_{i_1\ldots i_q})_{(i_1,\ldots,i_q)\in
R_1^{\prime}}$. Thus, $T$ is a vector space of dimension
$\begin{pmatrix} m
\\ q  \end{pmatrix}-\begin{pmatrix} m_1 \\
q \end{pmatrix}$.

Define $\Psi_1:L(\RR^m,\RR^m)\times T\to (\wedge^q\RR^m)^*$ by
$\Psi_1(\tau,F)=\tau^*F$ for $\tau\in L(\RR^m,\RR^m)$ and $F\in
T$. Since $(\tau_1^{-1})^*F_0\in T$, we have $F_0\in {\rm
im}(\Psi_1)$. Observe that for each $(i_1,\ldots,i_q)\in R$, the
$(i_1,\ldots,i_q)$ component of $\Psi_1(\tau,F)$, $(\tau^*F)_{i_1
\ldots i_q}$, is a polynomial of $\tau_{ji}$'s and $F_{j_1 \ldots
j_q}$'s, where $1\le i,j\le m$ and $(j_1,\ldots,j_q)\in
R_1^{\prime}$.
Identify $(\wedge^q\RR^m)^*$ with $\RR^{a_0}$ through standard basis $\{ e_1, \ldots, e_m\}$, where $a_0=\begin{pmatrix}m \\
q
\end{pmatrix}$. Note that since $q\ge 3$, if $m$ is sufficiently large, then
$m^2+\begin{pmatrix}
m \\ q  \end{pmatrix}-\begin{pmatrix} m_1 \\
q \end{pmatrix} < \begin{pmatrix} m \\ q
\end{pmatrix}$. Thus by Lemma~\ref{lem:algdep}, we can find a
non-zero polynomial $P_1$ such that
$$ {\rm im}(\Psi_1) \subset \{x \in \RR^{a_0}| P_1(x)=0 \} .$$

\underline{Case 2}. Suppose that $i_0=2,3$ or $4$. Let
$\tau_2:\RR^m\to \RR^m$ be the isomorphism sending $u_i$ to $e_i$
for $1\le i\le m$. Let $m_2$ denote $[\frac{m}{18}]$, then $m_2
\leq l$ and for $1\le i\le m_2$,
$$(\tau_2\varphi\tau_2^{-1})(e_i) = \lambda_i e_i +
\sum_{j=m_2+1}^m d_{ji} e_j,$$ for some $\lambda_i$ and
$d_{ji}$'s. Let $R_2=\{ (i_1, \ldots, i_q) | 1\leq i_1< \ldots<
i_q \leq m_2 \}$ and $R_2^{\prime}=R\backslash R_2$. Note that for
$(i_1,\ldots,i_q)\in R_2$, we have
$\lambda_{i_1}\ldots\lambda_{i_q}\ne 1$. (If $i_0=2$, then
$\lambda_{i_1}\ldots\lambda_{i_q}=-1$ since $q$ is odd; if
$i_0=3$, then $|\lambda_{i_1}\ldots\lambda_{i_q}|>1$; if $i_0=4$,
then $|\lambda_{i_1}\ldots\lambda_{i_q}|<1$.) Denote
$(\tau_2^{-1})^*F_0$ by $G$. Then for $(i_1,\ldots,i_q)\in R_2$,
$$F_0(u_{i_1},\ldots,u_{i_q})
=F_0(\varphi(u_{i_1}),\ldots,\varphi(u_{i_q})),$$ implies
$$G_{i_1\ldots
i_q}=G(e_{i_1},\ldots,e_{i_q})=G((\tau_2\varphi\tau_2^{-1})(e_{i_1}),\ldots,
(\tau_2\varphi\tau_2^{-1})(e_{i_q}))$$
$$=\lambda_{i_1}\ldots\lambda_{i_q}G_{i_1\ldots i_q}+\sum_{(j_1,
\ldots,j_q)\in R_2^{\prime}}c_{i_1 \ldots i_q, j_1 \dots j_q}
G_{j_1 \ldots j_q},$$ or equivalently,
\begin{equation}
G_{i_1\ldots i_q}= \frac{\sum_{(j_1, \ldots,j_q)\in
R_2^{\prime}}c_{i_1\ldots i_q, j_1 \ldots j_q} G_{j_1 \ldots
j_q}}{1-\lambda_{i_1}\ldots\lambda_{i_q}}, \label{eq:rational}
\end{equation}
where each coefficient $c_{i_1 \ldots i_q, j_1 \ldots j_q}$ is a
polynomial of $\lambda_i$'s and $d_{ji}$'s. Hence for any $(i_1,
\ldots, i_q) \in R_2$, $G_{i_1\ldots i_q}$ can be expressed as a
rational function of $\lambda_i$'s , $d_{ji}$'s and $G_{j_1 \ldots
j_q}$'s, where $1\le i\le m_2,m_2+1\le j\le m$ and
$(j_1,\ldots,j_q)\in R_2^{\prime}$.

Let $\Delta$ denote the open subset of $\RR^{m_2}$ consisting of
all $m_2$-tuples $(\lambda_1,\ldots,\lambda_{m_2})\in\RR^{m_2}$
such that for each $(i_1,\ldots,i_q)\in R_2$,
$\lambda_{i_1}\ldots\lambda_{i_q}\ne 1$. Let $M_{(m-m_2)\times
m_2}(\RR)$ denote the vector space consisting of all
$(m-m_2)\times m_2$ matrices with entries in $\RR$. Denote
$\begin{pmatrix}
m \\ q  \end{pmatrix}-\begin{pmatrix} m_2 \\
q \end{pmatrix}$ by $b_0$.

Define $\Psi_2:L(\RR^m,\RR^m)\times\Delta\times M_{(m-m_2)\times
m_2}(\RR)\times \RR^{b_0}\to (\wedge^q\RR^m)^*$ as follows: Let
$\tau\in
L(\RR^m,\RR^m),\lambda=(\lambda_1,\ldots,\lambda_{m_2})\in \Delta,
d=(d_{ji})_{m_2+1\le j\le m,1\le i\le m_2}\in M_{(m-m_2)\times
m_2}(\RR), f=(f_{j_1\ldots j_q})_{(j_1,\ldots, j_q)\in
R_2^{\prime}}\in \RR^{b_0}$. Then there is a unique $F\in
(\wedge^q\RR^m)^*$ such that for each $(j_1,\ldots,j_q)\in
R_2^{\prime}$, $F_{j_1\ldots j_q}=f_{j_1\ldots j_q}$ and for each
$(i_1,\ldots,i_q)\in R_2$, $F$ satisfies
equation~(\ref{eq:rational}) (changing $G$ to $F$ in that
equation). Define $\Psi_2(\tau,\lambda,d,f)$ to be $\tau^*F$. Note
that $F_0\in {\rm im}(\Psi_2)$. Observe that for each
$(i_1,\ldots,i_q)\in R$, the $(i_1,\ldots,i_q)$ component of
$\Psi_2(\tau,\lambda,d,f)$ is a rational function of
$\tau_{ji}$'s, $\lambda_i$'s, $d_{ji}$'s and $f_{j_1 \ldots
j_q}$'s.
Identify $(\wedge^q\RR^m)^*$ with $\RR^{a_0}$, where $a_0=\begin{pmatrix}m \\
q
\end{pmatrix}$. Note that since $q\ge 3$, if $m$ is sufficiently large, then
$m^2+m_2+(m-m_2)m_2+\begin{pmatrix}
m \\ q  \end{pmatrix}-\begin{pmatrix} m_2 \\
q \end{pmatrix} < \begin{pmatrix} m \\ q
\end{pmatrix}$. Thus by Lemma~\ref{lem:algdep}, we can find a
non-zero polynomial $P_2$ such that
$$ {\rm im}(\Psi_2) \subset \{x \in \RR^{a_0}| P_2(x)=0 \} .$$

For sufficiently large $m$, let $P=P_1P_2$ and $X_m= \{x \in
\RR^{a_0} | P(x)=0 \}$. $X_m$ is obviously closed. It cannot be
the whole $\RR^{a_0}$ since $P$ is a non-zero polynomial. Thus
$X_m$ is a proper closed subset of $(\wedge^q\RR^m)^*$. If $F\in
(\wedge^q\RR^m)^*$ is special, then $F\in {\rm im}(\Psi_1)\cup
{\rm im}(\Psi_2)\subset X_m$, as desired.

(b) We only need to substitute $(\wedge^q\RR^m)^*$ by
$(\vee^q\RR^m)^*$, substitute $\begin{pmatrix} m_i \\ q
\end{pmatrix}$ by $\begin{pmatrix} m_i+q-1 \\ q
\end{pmatrix}$ where $i= \varnothing, 1$ or $2$ and
change the restriction  $1 \le i_1 < \ldots < i_q \le m$ to $1\le
i_1\le \ldots\le i_q\le m$ in the definition of $R$, $R_1$ and
$R_2$. \hfill $\Box$

\section{Main Theorem}
\label{section:main}

In this section we prove the main theorem.

\begin{thm}
\label{thm:embed} If $n$ is a composite number and is not a power
of $2$, then there does not exist a single connected, closed
$(n+1)$-manifold $W$, such that any smooth, simply-connected,
closed $n$-manifold $M$ can be topologically flat embedded into
$W$.
\end{thm}

We only need to prove

\textbf{Theorem~\ref{thm:embed}*.} \textit{If $n$ is a composite
number and is not a power of $2$, then there does not exist a
single oriented, connected, closed $(n+1)$-manifold $W$, such that
any smooth, simply-connected, closed $n$-manifold $M$ can be
topologically flat embedded into $W$.}

\begin{defn} A \textit{graded algebra} over $\QQ$ is an algebra $A$ over $\QQ$
which as a vector space can be expressed as a direct sum
$$A=\bigoplus_{i=0}^{\infty}A_i,$$ and such that the bilinear
multiplication $\cup: A \times A \to A$ maps $A_r\times A_s$ into
$A_{r+s}$ for all nonnegative integers $r,s$. $A$ is
\textit{associative} if $(a\cup b)\cup c=a\cup (b\cup c)$ for all
$a,b,c\in A$, and is \textit{commutative} if for all $a\in
A_r,b\in A_s(r,s\ge 0)$, $a\cup b=(-1)^{rs}(b\cup a)$. We always
assume that $A$ is associative, $A_0=\QQ$ and $1\in\QQ =A_0$ is a
unit element of the algebra $A$.
\end{defn}

Note that the rational cohomology ring $H^*(X;\QQ)$ of a
path-connected topological space $X$ is a commutative graded
algebra over $\QQ$ (cf. [Sp], pp. 263-264).

\begin{defn} Let
$$A=\bigoplus_{i=0}^{\infty}A_i$$ be a commutative graded
algebra over $\QQ$. Assume that there exists a positive integer
$n$ such that $A_n\cong\QQ$, $A_i=0$ for $i>n$ and $A_i$ is finite
dimensional for $i<n$. Let $\varphi_i$ be the map from $A_i$ to
${\rm Hom}(A_{n-i}, A_n)$ defined by
$$\varphi_i(u): v \mapsto u\cup v, \ \forall  u\in A_i, \ v\in A_{n-i}.$$

We say that $A$ satisfies \textit{Poincar\'e duality} if
$\varphi_i$ is an isomorphism from $A_i$ to ${\rm Hom}(A_{n-i},
A_n)$ for any $0 \leq i \leq n$.
\end{defn}

\begin{prop}
\label{prop:extend} Suppose that $p$ is a positive integer, $q\ge
3$ is an odd integer. $F$ is in $(\wedge^qV)^*$ or $(\vee^qV)^*$
according to $p$ is odd or even, where $V$ is a finite dimensional
vector space over $\QQ$.

There exists a commutative graded algebra
$$A=\bigoplus_{i=0}^{pq} A_i$$ satisfying Poincar\'e duality such that

(i) $A_p=V$, $A_{pq}=\QQ$, $A_i=0$ if $i> pq$ and $A_i\not =0$
only if $i$ is a multiple of $p$,

(ii) for all $v_1,\ldots,v_q\in A_p=V$, $F(v_1, \ldots,  v_q)=
v_1\cup\ldots\cup v_q$.
\end{prop}

\begin{proof}
First assume that $p$ is odd. Let $q=2r+1$. We construct $A$
by setting
$$A_i=\left \{ \begin{array}{ll} \wedge^sV & {\rm if}\ i=sp\ {\rm
and}\ 0\le s\le r, \\
(\wedge^{2r+1-s}V)^* & {\rm if}\ i=sp\ {\rm and}\ r+1\le s\le
2r+1, \\ 0 & {\rm otherwise}
\end{array} \right.$$

Now suppose $u \in A_{sp}$, $v \in A_{tp}$ ($0\le s,t\le 2r+1$).
We define $u\cup v$, the product of $u$ and $v$, as follows. (We
assume that $s+t\le 2r+1$, since otherwise $u\cup v=0$.)

Case 1.   If $s+t\le r$, then $u\in A_{sp}=\wedge^sV, v \in
A_{tp}=\wedge^tV$, and we define $u\cup v\in
A_{(s+t)p}=\wedge^{s+t}V$ to be $u\wedge v$.

Case 2.   If $s\le r,t\le r$ and $s+t\ge r+1$, then we define
$u\cup v\in A_{(s+t)p}=(\wedge^{2r+1-s-t}V)^*$ by the equation
$(u\cup v)(w)=F(u\wedge v\wedge w)$ for $w\in \wedge^{2r+1-s-t}V$.

Case 3.   If $s\le r,t\ge r+1$ and $s+t\le 2r+1$, then $u\in
A_{sp}=\wedge^sV, v \in A_{tp}=(\wedge^{2r+1-t}V)^*$, and we
define $u\cup v\in A_{(s+t)p}=(\wedge^{2r+1-s-t}V)^*$ by the
equation $(u\cup v)(w)=v(w\wedge u)$ for $w\in
\wedge^{2r+1-s-t}V$.

Case 4. If $s\ge r+1,t\le r$ and $s+t\le 2r+1$, then $u\in
A_{sp}=(\wedge^{2r+1-s}V)^*, v \in A_{tp}=\wedge^tV$, and we
define $u\cup v\in A_{(s+t)p}=(\wedge^{2r+1-s-t}V)^*$ by the
equation $(u\cup v)(w)=u(v\wedge w)$ for $w\in
\wedge^{2r+1-s-t}V$.

Routine calculation shows that this product makes $A$ a
commutative graded algebra satisfying Poincar\'e duality and
conditions~(i) except that $A_{pq}=\QQ^*$. By the canonical
isomorphism $f: \QQ^* \to \QQ$ (i.e. for $w\in\QQ^*$, $f(w)=w(1)$)
, we can identify $A_{pq}$ with $Q$, then it is easy to see that
condition (ii) holds.

If $p$ is even, we can just change all $\wedge$ in the above
construction by $\vee$.\end{proof}

From [Su], Theorem~(13.2), we have

\begin{thm}
\label{thm:realize}
Any commutative graded algebra $A$ over $\QQ$
satisfying Poincar\'e duality and $A_1=0$ may be realized as the
rational cohomology ring of a simply-connnected, closed smooth
manifold of dimension $n$ with possibly one singular point when
$n=4k$. Moreover, when $n=4k$, if $A_{2k}=0$, then we can remove
the singularity in a realizing manifold.
\end{thm}

\textbf{Proof of Theorem~\ref{thm:embed}*.} Suppose such an
$(n+1)$-manifold $W$ exists. Since $n$ is a composite number and
is not a power of 2, we have a factorization $n=pq$, where $p$ is
an integer $\ge 2$ and $q$ is an odd integer $\ge 3$. By
Proposition~\ref{prop:special}, we can select an integer $m\ge
3\beta_{p+1}(W)$ such that $X_m$ can be constructed in
$(\wedge^q\RR^m)^*$ (when $p$ is odd) or $(\vee^q\RR^m)^*$ (when
$p$ is even). Since $X_m$ is proper closed, we can pick an
$F\not\in X_m$ such that $F(e_{i_1},\ldots,e_{i_q})\in \QQ$ for
any $1\le i_1,\ldots,i_q\le m$.

Taking $V=\QQ^m$, by Proposition~\ref{prop:extend}, there exists a
commutative graded algebra $A$ satisfying Poincar\'e duality and
conditions~(i) and (ii). $p\ge 2$ guarantees that $A_1=0$. By
Theorem~\ref{thm:realize}, there exists a simply-connected,
closed, smooth $n$-manifold $M$ such that the rational cohomology
ring $H^*(M;\QQ)$ is isomorphic to $A$. (When $n=4k$, $A_{2k}=0$
as $2k$ is not an integral multiple of $p$. Thus we can remove the
singularity in a realizing manifold.) Denote this isomorphism by
$h_\QQ$, and define $h_\RR$ to be $$h_\QQ \otimes {\rm id}: H^*(M;
\RR)= H^*(M; \QQ) \otimes_{\QQ} \RR \to A \otimes_{\QQ} \RR.$$
Clearly $h_\RR$ gives isomorphisms $H^p(M; \RR) \to \RR^m$ and
$H^n(M; \RR) \to \RR$. Condition (ii) in Proposition 4.2 can be
transformed to
\begin{equation} F(h_{\RR}(x_1),\ldots,
h_{\RR}(x_q))=h_{\RR}(x_1\cup\ldots\cup x_q). \label{eq:tran}
\end{equation}

Note that $\beta_p(M)=m\ge 3\beta_{p+1}(W)$, hence
$\frac{1}{2}(\beta_p(M)-\beta_{p+1}(W)) \ge \frac{m}{3}.$ Since
$M$ can be topologically flat embedded into $W$, by Proposition
2.1, there exists a subspace $V$ of $H^p(M;\RR)$ and a linear
transformation $\varphi: V \to H^p(M;\RR)$ with no fixed non-zero
vectors such that
$$\dim V \geq \frac{1}{2}(\beta_p(M)-\beta_{p+1}(W))\ge \frac{m}{3}$$
and for any $x_1, \ldots, x_q \in V$, \begin{equation} x_1 \cup
\ldots \cup x_q = \varphi(x_1) \cup \ldots \cup \varphi(x_q).
\label{eq:cuppre} \end{equation} Let $U=h_{\RR}(V)$. Define
$\psi:U\to \RR^m$ by $\psi(y)=(h_{\RR}\circ\varphi\circ
h_{\RR}^{-1})(y)$ for $y\in U$. Then $\dim U\ge \frac{m}{3}$ and
$\psi:U\to \RR^m$ is a linear map with no fixed non-zero vectors.
By equations~(\ref{eq:tran}) and ({\ref{eq:cuppre}), for any $y_1,
\ldots,y_q\in U$, we have
$$F(y_1,\ldots,y_q)=F(\psi(y_1),\ldots,\psi(y_q)).$$
Thus $F$ is special. But $F\not\in X_m$, a contradiction.\hfill
$\Box$

\section{Embedding of 4-manifolds into 5-manifolds}
\label{section:dim4case}

We first list theorems about the classification of
simply-connected $4$-manifolds according to their intersection
forms. More details can be found in [Ki] or [GS].

For any compact, connected, oriented $4$-manifold $M$, the cup
product
$$ \cup: H^2(M, \partial M) \times H^2(M, \partial M) \to H^4(M,
\partial M), (\alpha, \beta) \mapsto \alpha \cup \beta$$
gives a symmetric bilinear form $Q_M$ over $\ZZ$ on $H_2(M)$
through the isomorphisms $H^2(M, \partial M) \cong H_2(M)$ and
$H^4(M,
\partial M) \cong \ZZ $. Clearly, $Q_M(a,b)=0$ if $a$ or $b$ is a torsion element.
So $Q_M$ descends to a symmetric bilinear form
on $H_2(M)/\textrm{Torsion} \cong \ZZ^r$.

For a given integral symmetric bilinear form $Q$ on a finitely
generated free abelian group $A$, the \textit{rank},
\textit{signature} and \textit{parity} of $Q$ are defined as
follows: the rank $\textrm{rk}(Q)$ of $Q$ is the dimension of $A$.
Extend and diagonalize $Q$ over $A\otimes_\ZZ \RR$, the number of
positive entries and the number of negative entries are denoted by
$b_2^+$ and $b_2^-$ respectively, and the difference $b_2^+-b_2^-$
is the signature $\sigma(Q)$ of $Q$. Finally, $Q$ is even if
$Q(a,a)$ is even for any $a\in A$ and $Q$ is odd otherwise.

By choosing a basis for $A$, one can express $Q$ as a symmetric
matrix. The determinant of this matrix does not depend on the
choice of the basis, and thus we denote it by $\det Q$. If $\det
Q=\pm 1$, we say $Q$ is \textit{unimodular}. It follows from
Poincar\'e duality that $Q_M$ is unimodular if $M$ is closed. We
say $Q$ is \textit{indefinite} if both $b_2^+$ and $b_2^-$ are
positive.

\begin{thm}
\label{thm:indef} (cf. [GS], Theorem~1.2.21) Suppose $Q$ is an
integral, indefinite, unimodular symmetric bilinear form. If $Q$
is odd, then it is isomorphic to $b_2^+\langle+1\rangle \oplus
b_2^-\langle-1\rangle$; if $Q$ is even, then it is isomorphic to
$\frac{\sigma(Q)}{8} E_8 \oplus \frac{{\rm
rk}(Q)-|\sigma(Q)|}{2}H$, where $E_8$ is an even form with ${\rm
rk}(E_8) = \sigma(E_8)=8$ and $H=
\begin{bmatrix} 0 & 1 \\ 1 & 0 \end{bmatrix}$.
\end{thm}

This theorem classifies indefinite unimodular forms, while the
next theorem is a classification result for simply-connected
closed $4$-manifolds.

\begin{thm}
\label{thm:classify} (cf. [Fr]) For every integral, unimodular
symmetric bilinear form $Q$, there exists a simply-connected,
closed, topological $4$-manifold $M$ such that $Q_M \cong Q$.

If $Q$ is even, this manifold is unique up to homeomorphism. If
$Q$ is odd, there are exactly two different homeomorphism types of
such manifolds. At most one of these two homeomorphism types
carries a smooth structure. Consequently, simply-connected,
closed, smooth $4$-manifolds are determined up to homeomorphism by
their intersection forms.

In the odd case, the two manifolds are distinguished by their
Kirby-Siebenmann triangulation invariants ${\rm KS} \in \ZZ_2$ .
One has non-zero invariant and cannot be smooth, while the other
has zero invariant and is stably smoothable.
\end{thm}

For a compact, connected, topological $4$-manifold $M$, its
Kirby-Siebenmann invariant (cf. [KS]) ${\rm KS}(M)\in \ZZ_2$ is
the stable smoothing obstruction. This invariant is additive in
the following sense (cf. [FQ], Section~10.2B): If there is a
homeomorphism of boundaries $\partial M\cong \partial N$, then
${\rm KS}(M\cup_{\partial} N)={\rm KS}(M)+{\rm KS}(N)$. Combined
with the fact that $D^4$ is smoothable, we have

\begin{thm}
\label{thm:KS} Suppose $M_1$ and $M_2$ are two connected, closed,
topological $4$-manifolds, then
$${\rm KS}(M_1 \# M_2) = {\rm KS}(M_1) + {\rm KS}(M_2).$$
\end{thm}

Next, we give a construction that will be used many times later.

\begin{lem}
\label{lem:rembed}
For any oriented, closed, connected
$4$-manifold $M$, there exists an oriented, closed, connected
$5$-manifold $W$ such that for any positive integer $r$, $\#r M$
can be topologically flat embedded into $W$.
\end{lem}

\begin{proof} We first construct the double connected sum of $M$. Let
$D_1^4$ and $D_2^4$ be two disjoint embedded $4-$balls in $M$, and
let $X$ be
$$(M-{\rm int}(D_1^4)-{\rm int}(D_2^4)) \cup_\varphi (M-{\rm int}(D_1^4)-{\rm int}(D_2^4)),$$
where $\varphi$ is an orientation-reversing homeomorphism of
$\partial D_1^4 \coprod \partial D_2^4$.

\begin{center}
\resizebox{!}{4cm}{\includegraphics*[10mm, 10mm][170mm,
62mm]{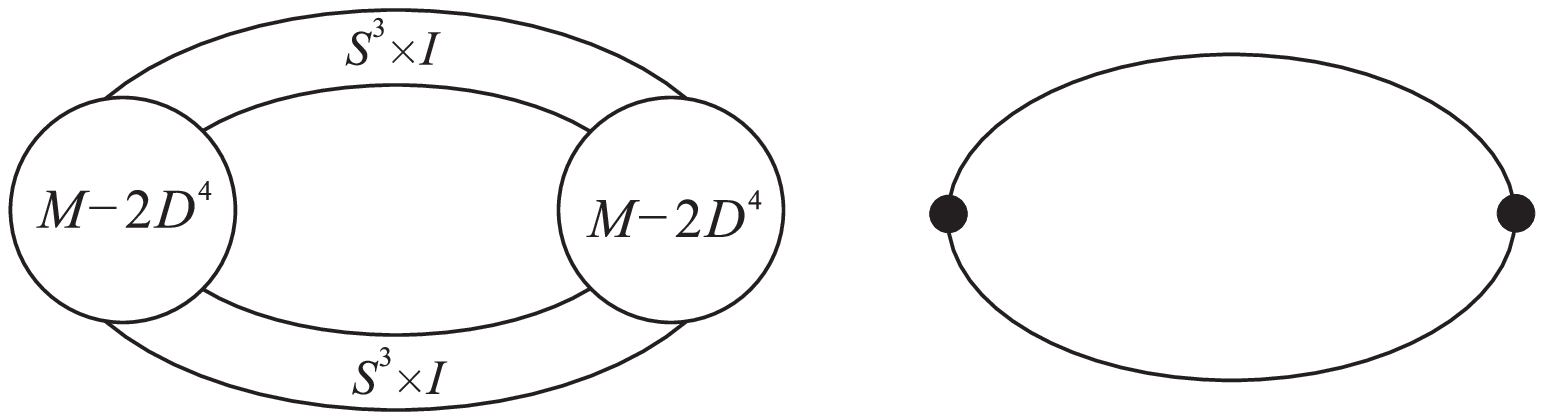} } Figure 1
\end{center}

Clearly, $X$ consists of 2 copies of $M - 2\ {\rm int}(D^4)$ and 2
copies of $S^3 \times I$. We denote these two kinds of pieces by
dots and arcs respectively as shown in Figure 1.

Next, we apply surgeries on $X \times S^1$ by cutting out two
disjoint copies of $(S^3 \times I ) \times I = S^3 \times D^2$ and
replacing each of them by $D^4 \times S^1$. The new manifold is
denoted by $W$.
\begin{center}
\resizebox{!}{6cm}{
\includegraphics*[10mm, 20mm][150mm,140mm]{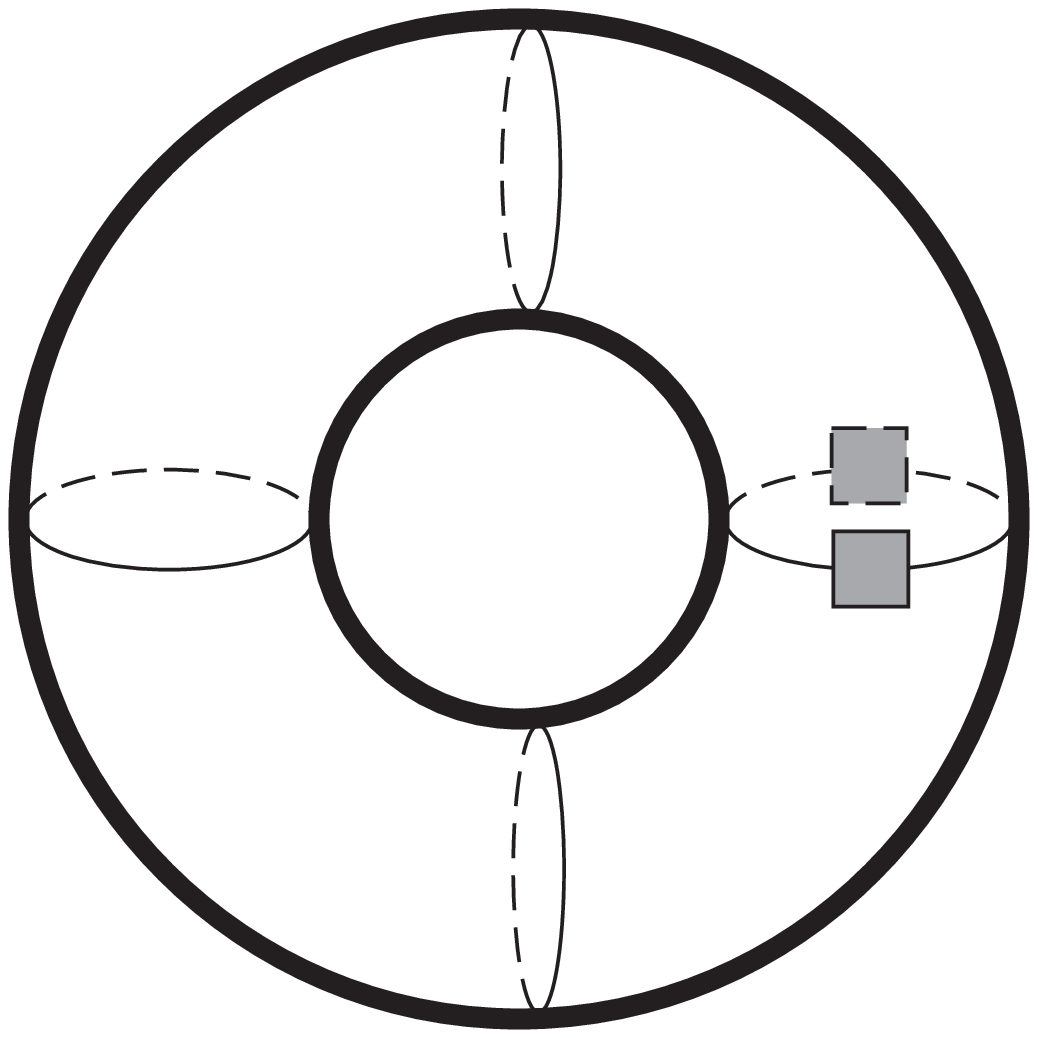} }

Figure 2
\end{center}

In Figure 2, $W$ is depicted as a torus, where the $X$ and the
$S^1$ factors correspond to meridian and longitude respectively.
In this picture, each point in the two thickened longitudes
represents $M - 2\ {\rm int}(D^4)$, any other point outside the
squares represents $S^3$, and the two shaded squares are the
surgery region.

We claim that for any positive integer $r$, $\#{r} M$ can be
embedded into $W$.

\begin{center}
\resizebox{15cm}{!}{ \includegraphics*[0mm,
0mm][230mm,100mm]{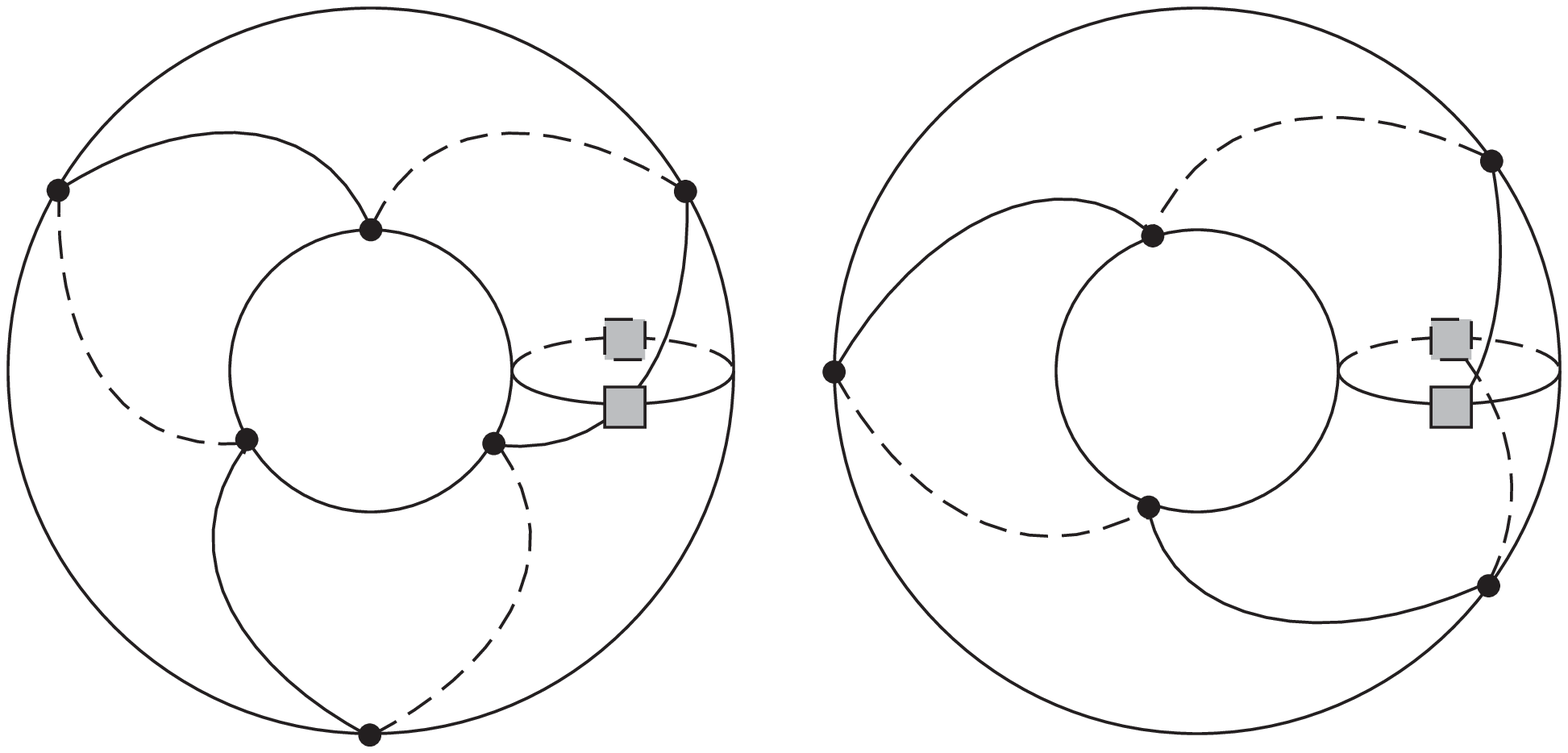} }  Figure 3
\end{center}

In Figure 3, a long curve is drawn with endpoints in the surgery
region and $r$ dots on it. (The left picture is for even $r$,
while the right picture is for odd $r$.) Each dot represents $ M -
2\ {\rm int}(D^4)$ and each arc in the curve that does not meet
the surgery region is $S^3 \times I$. For an arc that goes into
the surgery region, its intersection with the boundary of the
surgery region is $S^3$, which bounds a $D^4$ in $W$. Hence the
first and the last arc in the curve actually both denote $D^4$.
Now this curve gives an embedding of $\#{r} M$ in $W$.
\end{proof}

Now we can state and prove our results.

\begin{thm}
\label{thm:dim4case} (a) There exists an oriented, connected,
closed $5$-manifold $W$, such that for any simply-connected,
closed $4$-manifold $M$ with indefinite intersection form, $M$ can
be topologically flat embedded into $W$.

(b) Any simply-connected, compact topological $4$-manifold $M$
with $\partial M \not = \emptyset$ admits a topologically flat
embedding into $S^2 \tilde\times S^3$.
\end{thm}

\begin{proof} (a) Clearly, $\mathbb {CP} ^2$ and $\overline{\mathbb {CP}} ^2$
are smooth manifolds with intersection form $\langle+1\rangle$ and
$\langle-1\rangle$ respectively. And there are non-smoothable
manifolds $X_1$ and $X_{-1}$ that also have intersection form
$\langle+1\rangle$ and $\langle-1\rangle$ respectively. Then
$$ X_1\#{(p-1)}\mathbb {CP} ^2 \#{q}
\overline{\mathbb {CP}} ^2  \textrm{and} \ X_2\#{p}\mathbb {CP} ^2
\#{(q-1)} \overline{\mathbb {CP}} ^2$$ have intersection form
$p\langle+1\rangle \oplus q \langle-1\rangle$ and non-zero KS
invariant by Theorem~\ref{thm:KS} if the expressions are
meaningful. Now it follows from Theorems~\ref{thm:indef} and
\ref{thm:classify} that any simply-connected, oriented, closed
$4$-manifold $M$ with odd, indefinite intersection form is
homeomorphic to
$$k_1X_1 \# k_2X_{-1} \# k_3 \mathbb{CP} ^2 \# k_4
\overline{\mathbb {CP}} ^2,$$ for certain nonnegative integers
$k_i$.

Let $X_{E_8}$ and $X_{-E_8}$ be $4$-manifolds with intersection
form $E_8$ and $-E_8$, respectively. Clearly $S^2\times S^2$ has
intersection form $H$. Then it still follows from
Theorems~\ref{thm:indef} and \ref{thm:classify} that any
simply-connected, oriented, closed $4$-manifold $M$ with even,
indefinite intersection form is homeomorphic to
$$l_1X_{E_8} \# l_2X_{-E_8}\# l_3 (S^2 \times S^2) $$ for certain nonnegative integers
$l_i$.

Choose $M$ to be $X_1, X_{-1}, \mathbb{CP} ^2, \overline{\mathbb
{CP}} ^2, X_{E_8},X_{-E_8}$ and $S^2 \times S^2$ respectively and
construct corresponding $5$-manifolds $W_1, W_2 \ldots, W_7$ by
Lemma~\ref{lem:rembed}. It follows from the claim below that the
desired manifold $W$ can be chosen as $W_1\# W_2 \#\ldots \# W_7$.

\underline{Claim}. If there are topologically flat embeddings $M_1
\hookrightarrow W_1$ and $M_2 \hookrightarrow W_2$, then $M_1 \#
M_2$ can be topologically flat embedded into $W_1 \# W_2$.

\underline{Proof of the claim}. Choose $D^4_i$ in $M_i$. Since
$M_i\times I$ is a part of $W_i$, $D^4_i \times I$ is a $5$-ball
in $W_i$. Cut out these two $5$-balls from $W_i$ respectively and
glue their boundaries, then we have an embedding of $(M_1 \# M_2)
\times I$ into $W_1 \# W_2$.

(b) We glue $M$ and $-M$ together by the identity map between
$\partial M$ and $\partial(-M)$. The resulting manifold is usually
called the double of $M$ and denoted by $DM$. Clearly $DM$ is a
simply-connected, closed $4$-manifold with $\sigma(DM)=\sigma(M) +
\sigma(-M) =0$ (cf. [Ki], Chapter~II. Theorem~5.3) and ${\rm
KS}(DM)={\rm KS}(M)+{\rm KS}(-M)=0$. It follows from
Theorem~\ref{thm:indef} that the intersection form $Q$ of $DM$ is
either
$$ k\langle+1\rangle \oplus k\langle-1\rangle \  \textrm{or} \ \ k H,$$
depending on whether $Q$ is odd or
even.

Hence the intersection form of $DM$ is isomorphic to that of
either $\# k \mathbb{CP}^2 \#{k} \overline{\mathbb{CP}}^2$ or
$\#{k} S^2 \times S^2$. It is well known that
$$\mathbb{CP}^2 \# \overline{\mathbb{CP}}^2 \cong S^2
\tilde{\times} S^2, \ S^2  \tilde{\times} S^2 \# S^2
\tilde{\times} S^2 \cong S^2 \times S^2 \# S^2 \tilde{\times}
S^2.$$ Thus by ${\rm KS}(DM)=0$ and Theorem~\ref{thm:classify},
$DM$ is homeomorphic to either $\#{k} S^2 \times S^2$ or
$\#{(k-1)} S^2 \times S^2 \# S^2 \tilde{\times} S^2$.

$S^2 \tilde{\times} S^2$ can be embedded into $S^2 \tilde{\times}
S^3$ because the inclusion ${\rm SO}(3) \hookrightarrow {\rm
SO}(4)$ induces isomorphism $\pi_1({\rm SO(3)}) \cong \pi_1({\rm
SO}(4))  \cong \ZZ_2$. And $S^2 \times S^2$ can be embedded into
$\RR^5$ through the standard inclusion $S^2 \times D^3
\hookrightarrow \RR^5$. Therefore by the claim in the proof of
(a), both $\#{k} S^2 \times S^2$ and $\#{(k-1)} S^2 \times S^2 \#
S^2 \tilde{\times} S^2$ can be embedded into $S^2 \tilde{\times}
S^3$.
\end{proof}
\textbf{Acknowledgements.} The authors would like to thank
Jianzhong Pan for informing us the reference of Sullivan's work
[Su2] and for other interesting talks. The first two authors are
partially supported by grant no. 10201003 of NSFC and a grant of
MSTC.

\bigskip
\textsc{LMAM and School of Mathematical Sciences, Peking
University, Beijing 100871, P.R.China}

\textit{E-mail address:} \verb"dingfan@math.pku.edu.cn"

\bigskip
\textsc{LMAM and School of Mathematical Sciences, Peking
University, Beijing 100871, P.R.China}

\textit{E-mail address:} \verb"wangsc@math.pku.edu.cn"

\bigskip
\textsc{Department of Mathematics, University of California at
Berkeley, CA 94720, USA}

\textit{E-mail address:} \verb"jgyao@math.berkeley.edu"


\begin{thebibliography}{123}

\bibitem[Di]{Di} {\bf F. Ding}, {\it Smooth structures on some open $4$-manifolds}, Topology
36 (1997), no. 1, 203--207.


\bibitem[EL]{EL} {\bf A. L. Edmonds, C. Livingston}, {\it Embedding punctured lens
spaces in four-manifolds}, Comment. Math. Helv. 71 (1996), no. 2,
169--191.

\bibitem[Fa1]{Fa1}{\bf F. Fang}, {\it Embedding $3$-manifolds and smooth structures of
$4$-manifolds}, Topology Appl. 76 (1997), no. 3, 249--259.

\bibitem[Fa2]{Fa2}{\bf F. Fang}, {\it Smooth structures on
$\Sigma\times {\bf R}$}, Topology Appl. 99 (1999), no. 1,
123--131.

\bibitem[Fr]{Fr} {\bf M. H. Freedman}, {\it The topology of
four-dimensional manifolds}, J. Differential Geom. 17 (1982), no.
3, 357--453.

\bibitem[FQ]{FQ} {\bf M. H. Freedman, F. Quinn}, {\it Topology of
4-manifolds}, Princeton University Press, Princeton, NJ, 1990.

\bibitem[GS]{GS} {\bf R. E. Gompf, A. I. Stipsicz}, {\it $4$-manifolds
and Kirby calculus}, Graduate Studies in Math. 20, Amer. Math.
Soc., Providence, RI, 1999.

\bibitem[Hi]{Hi}{\bf M. W. Hirsch}, {\it The imbedding of bounding manifolds in euclidean
space}, Ann. of Math. (2) 74 (1961), 494--497.

\bibitem[HJ]{HJ} {\bf R. A. Horn, C. R. Johnson}, {\it Matrix
analysis}, Cambridge University Press, Cambridge, 1990.


\bibitem[Ka]{ka} {\bf A. Kawauchi}, {\it The imbedding problem of
$3$-manifolds into $4$-manifolds}, Osaka J. Math. 25 (1988), no.
1, 171--183.

\bibitem[Ki]{ki} {\bf R. C. Kirby}, {\it The topology of
$4$-manifolds}, Lecture Notes in Math. 1374, Springer-Verlag,
Berlin, 1989.

\bibitem[KS]{KS} {\bf R. C. Kirby, L. C. Siebenmann}, {\it Foundational
essays on triangulations and smoothings of topological manifolds,
}, Ann. of Math. Studies 88, Princeton University Press,
Princeton, NJ, 1977.

\bibitem[La]{la} {\bf S. Lang}, {\it Algebra}, Graduate Texts in Math. 211,
Springer-Verlag, New York, 2002.

\bibitem[Ro]{Ro}{\bf V. A. Rohlin}, {\it The embedding of non-orientable three-manifolds into
five-dimensional Euclidean space}, Soviet Math. Dokl. 6 (1965),
153--156.

\bibitem[Sh]{Sh} {\bf R. Shaw}, {\it Linear algebra and group
representations, Vol. II}, Academic Press, London-New York, 1983.

\bibitem[Shi]{Shi} {\bf T. Shiomi}, {\it On imbedding
$3$-manifolds into $4$-manifolds}, Osaka J. Math. 28 (1991), no.
3, 649--661.

\bibitem[Sp]{Sp} {\bf E. H. Spanier}, {\it Algebraic topology},
Springer-Verlag, New York, 1966.


\bibitem[Su]{Su} {\bf D. Sullivan}, {\it Infinitesimal computations in topology},
Inst. Hautes \'Etudes Sci. Publ. Math. 47 (1977), 269--331 (1978).

\bibitem[Wa]{Wa} {\bf C. T. C. Wall}, {\it All $3$-manifolds imbed in $5$-space}, Bull. Amer.
Math. Soc. 71 (1965), 564--567.

\bibitem[WZ]{WZ} {\bf S. Wang, Q. Zhou}, {\it How to embed $3$-manifolds into
$5$-space}, Adv. in Math. (China) 24 (1995), no. 4, 309--312.

\end{thebibliography}
\end{document}